\documentstyle[amssymb, 12pt]{article}

\begin{document}
\font\frak=eufm10 scaled\magstep1
\font\fak=eufm10 scaled\magstep2
\font\fk=eufm10 scaled\magstep3
\font\scriptfrak=eufm10
\font\tenfrak=eufm10

\newtheorem{theorem}{Theorem}
\newtheorem{corollary}{Corollary}
\newtheorem{proposition}{Proposition}
\newtheorem{definition}{Definition}
\newtheorem{lemma}{Lemma}
\font\frak=eufm10 scaled\magstep1
\newenvironment{pf}{{\noindent{\it Proof. }}}{\ $\Box$\medskip}


\mathchardef\za="710B  
\mathchardef\zb="710C  
\mathchardef\zg="710D  
\mathchardef\zd="710E  
\mathchardef\zve="710F 
\mathchardef\zz="7110  
\mathchardef\zh="7111  
\mathchardef\zvy="7112 
\mathchardef\zi="7113  
\mathchardef\zk="7114  
\mathchardef\zl="7115  
\mathchardef\zm="7116  
\mathchardef\zn="7117  
\mathchardef\zx="7118  
\mathchardef\zp="7119  
\mathchardef\zr="711A  
\mathchardef\zs="711B  
\mathchardef\zt="711C  
\mathchardef\zu="711D  
\mathchardef\zvf="711E 
\mathchardef\zq="711F  
\mathchardef\zc="7120  
\mathchardef\zw="7121  
\mathchardef\ze="7122  
\mathchardef\zy="7123  
\mathchardef\zf="7124  
\mathchardef\zvr="7125 
\mathchardef\zvs="7126 
\mathchardef\zf="7127  
\mathchardef\zG="7000  
\mathchardef\zD="7001  
\mathchardef\zY="7002  
\mathchardef\zL="7003  
\mathchardef\zX="7004  
\mathchardef\zP="7005  
\mathchardef\zS="7006  
\mathchardef\zU="7007  
\mathchardef\zF="7008  
\mathchardef\zW="700A  

\newcommand{\be}{\begin{equation}}
\newcommand{\ee}{\end{equation}}
\newcommand{\ra}{\rightarrow}
\newcommand{\lra}{\longrightarrow}
\newcommand{\bea}{\begin{eqnarray}}
\newcommand{\eea}{\end{eqnarray}}
\newcommand{\beas}{\begin{eqnarray*}}
\newcommand{\eeas}{\end{eqnarray*}}
\newcommand{\Z}{{\Bbb Z}}
\newcommand{\R}{{\Bbb R}}
\newcommand{\C}{{\Bbb C}}
\newcommand{\1}{{\bold 1}}
\newcommand{\SL}{SL(2,\C)}
\newcommand{\Sl}{sl(2,\C)}
\newcommand{\SU}{SU(2)}
\newcommand{\su}{su(2)}
\newcommand{\G}{{\goth g}}
\newcommand{\D}{{\rm d}}
\newcommand{\de}{\,{\stackrel{\rm def}{=}}\,}
\newcommand{\we}{\wedge}
\newcommand{\nn}{\nonumber}
\newcommand{\ot}{\otimes}
\newcommand{\s}{{\textstyle *}}
\newcommand{\ts}{T^\s}
\newcommand{\da}{\dagger}
\newcommand{\pa}{\partial}
\newcommand{\ti}{\times}
\newcommand{\A}{{\cal A}}
\newcommand{\Li}{{\cal L}}
\newcommand{\ka}{{\Bbb K}}
\newcommand{\find}{\mid}

\title{Non-antisymmetric versions\\ of  Nambu-Poisson
and Lie algebroid brackets}

\author{
Janusz Grabowski\thanks{Supported by KBN, grant No. 2 P03A 031 17.}\\
Institute of Mathematics, Warsaw University\\
ul. Banacha 2, 02-097 Warszawa, Poland. \\
and\\
Mathematical Institute, Polish Academy of Sciences\\
ul. \'Sniadeckich 8, P. O. Box 137, 00-905 Warszawa, Poland\\
{\it e-mail:} jagrab@mimuw.edu.pl
\and
Giuseppe Marmo\thanks{Supported by PRIN SINTESI}\\
Dipartimento di Scienze Fisiche,
Universit\`a Federico II di Napoli\\
and\\
INFN, Sezione di Napoli\\
Complesso Universitario di Monte Sant'Angelo\\
Via Cintia, 80126 Napoli, Italy\\
{\it e-mail:} marmo@na.infn.it}
\maketitle
\begin{abstract}
We show that we can skip the skew-symmetry assumption in
the definition of Nambu-Poisson brackets. In other words, a
$n$-ary bracket on the algebra of smooth functions which
satisfies the Leibniz rule and a $n$-ary version of the
Jacobi identity must be skew-symmetric. A similar result
holds for a non-antisymmetric version of Lie algebroids.
\end{abstract}

\section{Introduction}
Two main directions have been suggested for the generalization of the
notion of a Lie algebra. First, Filippov developed a proposal for
brackets with more than two arguments, i.e., $n$-ary brackets. In
\cite{Fi}
he proposed a definition of such structures (which we shall call {\it
Filippov algebras}) with a version of the Jacobi identity
for $n$-arguments (we shall call it {\it Filippov identity}, shortly FI):
\be\label{FI}
\{f_1,\dots,f_{n-1},\{ g_1,\dots,g_n\}\}=\sum_{k=1}^n
\{g_1,\dots,\{f_1,\dots,f_{n-1},g_k\},\dots,g_n\}.
\ee
Note that in the binary case ($n=2$), the Filippov identity coincides with
the Jacobi identity.
Independently, Nambu \cite{Nam}, looking for generalized formulations
of Hamiltonian Mechanics, found $n$-ary analogs of Poisson
brackets and then Takhtajan \cite{Ta} rediscovered the Filippov
identity (and called it Fundamental Identity) for them.
The {\it Filippov brackets} are assumed to be $n$-linear and
skew-symmetric and {\it Nambu-Poisson brackets}, defined on algebras of
smooth functions, satisfy additionally the {\it Leibniz rule}:
\be
\{f_1f_1',\dots,f_n\}=f_1\{f_1',\dots,\dots,f_n\}+
\{f_1,\dots,\dots,f_n\}f_1'.
\ee
On the other hand, Loday (cf. \cite{Lo1}), while studying relations between
Hochschild and cyclic homology in the context of searching for
obstructions to the periodicity of algebraic K-theory, discovered that
one can skip the skew-symmetry assumption in the definition of a Lie algebra,
still having a possibility to define appropriate (co)homology
(see \cite{Lo1,LP} and \cite{Lo}, Chapter 10.6). His Jacobi
identity for such structures was formally the same as the classical Jacobi
identity in the form of (\ref{FI}) for $n=2$:
\be\label{JI}
\{ f,\{ g,h\}\}=\{\{ f,g\},h\}+\{ g,\{ f,h\}\}.
\ee
This time, however, this is no longer equivalent to
\be\label{JI1}
\{\{ f,g\},h\}=\{\{ f,h\},g\}+\{ f,\{g,h\}\},
\ee
since we have no skew-symmetry. Loday called such structures {\it
Leibniz algebras} but, since we have already associated the name of
Leibniz with the Leibniz
identity, we shall call them {\it Loday algebras}. This is in accordance
with the terminology of \cite{KS}, where
analogous structures in the graded case are defined. Of course, there is
no particular reason not to define Loday algebras by means of (\ref{JI1})
instead of (\ref{JI}) (and in
fact, it was the original Loday definition), but both categories are
equivalent via transposition of arguments. Similarly, for associative
algebras we can obtain associated algebras by transposing arguments, but
in this case we still get associative algebras. It is interesting that
Nambu-Poisson brackets lead to some Loday algebras and hence to
the corresponding (co)homology (see \cite{DT}).
\par
It is now clear that we can combine both generalizations and define {\it
Filippov-Loday} algebras as those which are  equipped with  $n$-ary
brackets, not
skew-symmetric in general, but satisfying the Filippov identity. We can
also define a Loday version of Nambu-Poisson algebras or rings (we
shall call them {\it Nambu-Poisson-Loday}, or simply {\it Nambu-Loday},
algebras (or rings)), assuming additionally that a Filippov-Loday
structure is defined on a commutative associative algebra
(resp. ring) and satisfies the Leibniz rule (with respect to all
arguments separately, since we have no skew-symmetry).
\par
In this short note we first deal with the problem of finding examples of
new, i.e., non antisymmetric, Nambu-Poisson-Loday brackets. The result is,
to some extend, unexpected. We show that for a wide variety of associative
commutative algebras, including algebras of smooth functions, we
get nothing more than what we already know, since Nambu-Loday algebras
have to be skew-symmetric.
In particular, we can skip the skew-symmetry axiom in the standard
definition of Poisson bracket.

We obtain a similar negative result for a Loday-type generalization
of Lie algebroids: they are locally, in principle, skew-symmetric, or
they are bundles of Loday algebras.

\section{Main Theorem}
\noindent
{\bf Definition.} Let $\A$ be an associative commutative algebra.
Let $\{\cdot,\dots,\cdot\}$ be a {\it $n$-ary bracket} on $\A$, i.e.,
an operation with $n$-arguments.
\be
\A\ti\cdots\ti\A\ni(f_1,\dots,f_n)\mapsto\{ f_1,\dots,f_n\}\in\A
\ee
which is linear with respect to all arguments:
\be
\{f_1,\dots,\za f_i+\zb f_i',\dots,f_n\}=\za\{f_1,\dots,f_i,\dots,f_n\}+
\zb\{f_1,\dots,f_i',\dots,f_n\}.
\ee
We shall call such a bracket a {\it Nambu-Loday} bracket, if it satisfies
the following two conditions:
\par
(i) the Leibniz rule (with respect to each
argument):
\be
\{f_1,\dots,f_if_i',\dots,f_n\}=f_i\{f_1,\dots,f_i',\dots,f_n\}+
\{f_1,\dots,f_i,\dots,f_n\}f_i'\ ,
\ee
\par
for all $i=1,\dots,n$, and
\par
(ii) the Filippov identity:
\be\label{fi}
\{f_1,\dots,f_{n-1},\{ g_1,\dots,g_n\}\}=\sum_{k=1}^n
\{g_1,\dots,\{f_1,\dots,f_{n-1},g_k\},\dots,g_n\}.
\ee
The commutative algebra $\A$, equipped with a Nambu-Loday bracket, will be
called {\it Nambu-Loday algebra}. For Nambu-Loday algebras we have no direct
inductive characterization as that for Nambu-Poisson and Nambu-Jacobi
brackets \cite{GM}, since the property saying that fixing an argument
we get a bracket satisfying FI, but of one argument less, is based on
skew-symmetry. However, we can prove the following.
\begin{theorem} If $\A$ is an associative commutative algebra over a field
of characteristic 0 and $\A$ contains no
nilpotents, then every Nambu-Loday bracket on $\A$ is skew-symmetric.
\end{theorem}
\begin{pf} Let us assume that we have fixed a Nambu-Loday bracket
on $\A$.
First, observe that the skew-symmetry property is equivalent to the fact
that the bracket vanishes, if only two arguments are the same. Explicitly,
if
\be\label{10}
\{ f_1,\dots,f_n\}=0 {\rm\ for\ all\ } f_1,\dots,f_n\in\A {\rm \ with\ }
f_i=f_j=h {\rm \ for\ some\ } i\ne j,
\ee
then writing $h=x+y$ and using (\ref{10}) for $h=x$ and $h=y$, we get
\bea\label{2}
&&\{ f_i,\dots,f_{i-1},x,f_{i+1},\dots,f_{j-1},y,f_{j+1},\dots,f_n\}=\\
&&-\{ f_i,\dots,f_{i-1},y,f_{i+1},\dots,f_{j-1},x,f_{j+1},\dots,f_n\}. \nn
\eea
Second, since we can get the skew-symmetry (\ref{2}) with respect to the
transposition $(i,j)$ composing transpositions $(i,n),(j,n)$ and $(i,n)$
again, it is sufficient to prove (\ref{2}) (or (\ref{10})) for $j=n$.
\par
Fix $i=1,\dots,n-1$. Replacing $f_i$ in (\ref{fi}) by $f_i^2/2$, we
get, due to the Leibniz rule,
\be\label{1}
f_i\{f_1,\dots,f_{n-1},\{ g_1,\dots,g_n\}\}=\sum_{k=1}^n
\{g_1,\dots,f_i\{f_1,\dots,f_{n-1},g_k\},\dots,g_n\}.
\ee
Subtracting from (\ref{1}) the Filippov identity (\ref{fi}) multiplied by
$f_i$, we get
\be\label{L}
\sum_{k=1}^n\left(\{ f_1,\dots,f_{n-1},g_k\}\{ g_1,\dots,g_{k-1},f_i,
g_{k+1},\dots,g_n\}\right)=0,
\ee
which holds for all $i=1,\dots,n-1$
\par
Now, for $n-1\ge m\ge 1$, we shall show inductively the following:
\be
(S_m)\quad {\sl \ If\ } m {\sl \ elements\ of\ } f_1,\dots,f_{n-1}\in\A
{\sl \ equal\ } h, {\sl \ then \ } \{ f_1,\dots,f_{n-1},h\}=0.
\ee
Of course, $(S_1)$ tells us just that the bracket is skew-symmetric
with respect to all transpositions $(i,n)$, so it is totally
skew-symmetric, according to the previous remarks.

We start with $m=n-1$. Putting in (\ref{L}) all $f$'s and $g$'s equal
to $h$, we get $n\{ h,\dots,h\}^2=0$, which gives us $\{ h,\dots,
h\}=0$, since there are no nilpotents in $\A$ , so the induction
starts.  To prove the inductive step, assume ($S_m$) for some $n-1\ge
m>1$. We shall show ($S_{m-1}$). Take $f_1,\dots,f_n\in\A$ such that
$f_j=h$ for $j$ from a subset $I$ of $\{ 1,\dots, n-1\}$ with $(m-1)$
elements. Put $f_k=g_k$, $k=1,\dots,n-1$, and $g_n=h$.
For a fixed $i=1,\dots,n-1$, we have
(i) if $k\notin I$, then $\{ g_1,\dots,g_{k-1},f_i,g_{k+1},\dots,g_n\}=0$
by the inductive assumption, and

(ii) if $k\in I$, then
$$\{ f_1,\dots,f_{n-1},g_k\}\{ g_1,\dots,g_{k-1},f_i,
g_{k+1},\dots,g_n\}=\{ f_1,\dots,f_{n-1},h\}^2.
$$
This implies that (\ref{L}) reads in this case
\be
m\{ f_1,\dots,f_{n-1},h\}^2=0,
\ee
which gives
\be
\{ f_1,\dots,f_{n-1},h\}=0,
\ee
for any $f_1,\dots,f_{n-1}\in\A$ such that $(m-1)$ of them equal $h$.
\end{pf}

\medskip\noindent
{\bf Remark.} The assumption that there are no nilpotents in $\A$
is essential. To see this, consider the commutative associative
algebra $\A$ over ${\frak k}$ freely
generated by $x,y$, with the constraint $x^2=0$, i.e., $\A={\frak
k}[x,y]/\langle x^2\rangle$. It is easy to see that the bracket
\be
\{ f,g\}=x\frac{\pa f}{\pa y}\frac{\pa g}{\pa y}
\ee
satisfies the Leibniz rule and the Jacobi identity, but it is symmetric.

\section{Generalized Lie algebroids}

Every $n$-ary bracket on the algebra $C^\infty(M)$ of smooth functions
on a manifold $M$, which satisfies the Leibniz rule, is associated with
a $n$-contravariant tensor $\zL$ according to
\be\label{lin}
\{ f_1,\dots,f_n\}=\langle\zL,\D f_1\we\cdots\we\D f_n\rangle.
\ee
The vector fields
\be
\zL_{(f_1,\dots,f_{n-1})}=\{ f_1,\dots,f_{n-1},\cdot\}=i_{\D
f_1\we\cdots\we\D f_{n-1}}\zL
\ee
we can call (left) {\it Hamiltonian vector fields} of $\zL$.
It is easy to see that the Filippov identity for the $n$-bracket is in
this case equivalent to the fact that the Hamiltonian vector fields
preserve the tensor $\zL$, i.e.,
\be
\Li_{\zL_{(f_1,\dots,f_{n-1})}}\zL=0,
\ee
where $\Li$ stands for the Lie derivative.
Theorem 1 can be formulated in this case as follows.
\begin{corollary} If a $n$-contravariant tensor field is preserved by its
Hamiltonian vector fields, then it is skew-symmetric.
\end{corollary}
It is well known that with a  $n$-ary bracket on a finite-dimensional
vector space $V$ (over $\R$) we canonically associate a linear
contravariant $n$-tensor $\zL$ on the dual space $V^*$ such that
(\ref{lin}) is
satisfied for linear functions on $V^*$ (thus elements of $V$).
Explicitly, if $(x_1,\dots,x_k)$ is a basis of $V$ (thus a coordinate
system of $V^*$), then
\be
\zL=\sum_{i_1,\dots,i_n=1}^k\{ x_{i_1},\dots,x_{i_n}\}\pa_{i_1}\otimes
\cdots\otimes\pa_{i_n}.
\ee
Lie algebras correspond in this way to linear Poisson tensors. This can be
generalized to vector bundles as follows.
By linear functions on a vector bundle $E$ over a manifold
$M$ we understand the functions we get from sections of the dual bundle
$E^*$ by contraction, i.e., the linear function $\zi_X$ associated with a
section $X$ of $E^*$ is given by $\zi_X(\za_p)=\langle
X(p),\za(p)\rangle$, where $\za(p)\in E_p$ for $p\in M$. We say that a
$n$-tensor $\zL$ on $E^*$ is linear if linear functions on $E^*$ are
closed with respect to the $n$-ary bracket $\{\cdot,\dots,\cdot\}_\zL$
generated by $\zL$. Hence, we can define a $n$-ary operation
$[\cdot,\dots,\cdot]_\zL$ on sections of the bundle $E$ by
\be
\zi_{[X_1,\dots,X_n]_\zL}=\{\zi_{X_1},\dots,\zi_{X_n}\}_\zL.
\ee
In \cite{GU1,GU2} this idea was used to define general (binary) algebroid
structures, and in \cite{GM1} to define $n$-ary Lie algebroids.
Let us concentrate now on the binary case and let us recall from
\cite{GU2} the following definition.
\begin{definition} Let $M$ be a manifold. An {\rm  algebroid  on
$M$} is a vector bundle $\zt:E\ra M$, together with a  bracket
$[\cdot,\cdot]:\A\ti \A\ra\A$ on the module $\A=\zG E$ of
global sections of $E$, and two vector bundle morphism $a_l,a_r:E\ra TM$,
over the  identity on $M$, from $E$ to the tangent bundle $TM$, called the
{\rm anchors}  of the Lie algebroid (left and right), such that
\be
[fX,gY]=fg[X,Y]+fa_l(X)(g)Y-ga_r(Y)(f)X,
\ee
for all  $X,Y\in\A$  and  all $f,g\in\C(M)$.
\end{definition}
It is clear that any finite-dimensional algebra structure can be viewed as
an algebroid structure on a bundle over a single point.
Note that in the case when the algebroid bracket is a Lie bracket, we have
$a_l=a_r=a$ and
$a([X,Y])=[a(X),a(Y)]$ for all $X,Y\in\A$.
Such structures are called Lie algebroids. They
were introduced by Pradines \cite{Pr} as infinitesimal objects for
differentiable groupoids, but one can find similar notions proposed by
several authors in increasing number of papers (which proves their
importance and naturalness). For basic properties and the literature on the
subject we refer to the survey article by Mackenzie \cite{Ma}.

\begin{theorem}(\cite{GU1}) There is a one-one correspondence
between linear 2-contravariant tensors $\zL$ on the dual bundle $E^*$ and
algebroid brackets $[\cdot,\cdot]_\zL$ on $E$.
\end{theorem}
Note that, equivalently, we can think of algebroid structures on
the vector bundle $E$ as morphisms of double vector bundles $\ze:T^*E
\ra TE^*$ (cf. \cite{GU2}).

We can speak about {\it Loday algebroids} when we impose the Jacobi
identity (\ref{JI}) but we skip the skew-symmetry assumption.
One can think that imposing the Jacobi identity for an algebroid, we get
the Jacobi identity for the bracket $\{\cdot,\cdot\}_\zL$ of
functions defined by the
corresponding tensor $\zL$ on $C^\infty(E^*)$ and, in view of Theorem 1,
that this implies that $\zL$ is a Poisson tensor, so our algebroid is a
Lie algebroid. This reasoning, however, is wrong, since the Jacobi
identity on sections of $E$ forces the Jacobi identity for the bracket
$\{\cdot,\cdot\}_\zL$ only for linear functions. Such tensors may be
non-skew-symmetric, i.e., clearly, Loday algebras do exist.
A simple example is the following.

\medskip\noindent
{\bf Example 1.} Consider the 2-tensor on $\R^3$ given by
\be
\zL=x_2\pa_1\otimes\pa_1+x_3\pa_1\otimes\pa_3-x_3\pa_3\otimes\pa_1.
\ee
It is easy to see that the Hamiltonian vector fields of linear functions
preserve $\zL$, so we have the Jacobi identity for the associated bracket:
\be
[x_1,x_1]=x_2, \quad [x_1,x_3]=-[x_3,x_1]=x_3,
\ee
where we assume the missing brackets to be zero.
This example is also an example of a Loday algebroid over a single
point, but we can obtain a Loday algebroid over $M$ just tensoring the
above algebra with $C^\infty(M)$.

\medskip\noindent
The anchors of the Loday algebroids from the above example are trivial. We
shall show that this is not incidental and Loday algebroids can be reduced
to Lie algebroids and bundles of Loday algebras.
\begin{theorem}
For any Loday algebroid bracket the left anchor is the same as the right
anchor and the bracket is skew-symmetric at points where they do not
vanish.
\end{theorem}
\begin{pf} Let $[\cdot,\cdot]$ be a Loday algebroid bracket on the
space $\A$ of sections of a vector bundle $E$ over $M$. The Jacobi
identity implies immediately
\be\label{21}
[[X,X],Y]=0,
\ee
for all $X,Y\in\A$. Putting $X:=fX$ in (\ref{21}), we get
\bea\label{22}
&&f(a_l(X)(f)-a_r(X)(f))[X,Y]-2fa_r(Y)(f)[X,X]-\\
&&fa_r(Y)(a_l(X)(f)-a_r(X)(f))X-(a_l(X)(f)-a_r(X)(f))a_r(Y)(f)X=0,\nn
\eea
for all $X,Y\in\A$ and all $f\in C^\infty(M)$.
Suppose that at $p\in M$ the right anchor does not vanish, i.e., there
are $Y\in\A$ and $f\in C^\infty (M)$ such that $a_r(Y)(f)(p)\ne 0$. We
can additionally assume that $f(p)=0$ and then (\ref{22}) implies that
$(a_l(X)(f)-a_l(X)(f))(p)=0$ for all $X\in\A$. Hence the vector
$(a_l(X)-a_r(X))(p)$ annihilates any covector from $T^*_pM$ not
annihilated by $a_r(Y)(p)$, thus it is zero. But if $a_r$ does not vanish
at $p$, then it does not vanish in a neighborhood of $p$, so
$a_l(X)=a_r(X)$ in a neighborhood of $p$ and (\ref{22}) implies now
that in this neighborhood
\be
fa_r(Y)(f)[X,X]=0,
\ee
for all $X,Y\in\A$ and $f\in C^\infty (M)$. Since $a_r$ is nontrivial in
this neighborhood, this in turn implies $[X,X]=0$, i.e., the bracket is
skew-symmetric. In particular, the left anchor equals the right one.

Assume now that the right anchor vanishes at $p\in M$. By (\ref{22}) we
obtain now
\be
f(p)a_l(X)(f)(p)[X,Y](p)=0,
\ee
for all $X,Y\in\A$ and $f\in C^\infty (M)$, so
\be\label{23}
a_l(X)(f)(p)[X,Y](p)=0.
\ee
Replacing $X$ in (\ref{23}) by $X+Z$, we get
\be
a_l(X)(f)(p)[Z,Y](p)+a_l(Z)(f)[X,Y](p)=0.
\ee
Multiplying the above equation by $a_l(X)(f)(p)$ and taking into account
(\ref{23}), we get
\be
(a_l(X)(f)(p))^2[Z,Y](p)=0,
\ee
for all $X,Y,Z\in\A$ and $f\in C^\infty (M)$ which clearly implies that
the left anchor vanishes at $p$, since, if the bracket is trivial at $p$,
then both anchors are trivial at $p$.
Hence, the right anchor is the same as the left anchor and the bracket is
skew-symmetric at points where they do not vanish.
\end{pf}
\section{Conclusions}
Poisson and Lie algebroid brackets are ones of the most fundamental
algebraic structures in Classical and Quantum Physics.
We have composed the two ways of generalizing Poisson bracket:
the Nambu's idea of $n$-ary bracket and the Loday's observation that
skipping the skew-symmetry assumption in the definition of a Lie
algebra we still have a (co)homology theory. What we get is that no
new structures appear in this way, since the Leibniz rule and the
Filippov identity imply the skew-symmetry. A similar phenomena we find
out when looking for a non-skew version of a Lie algebroid.
This shows that skew-symmetry is in fact forced by other properties of
these important algebraic structures.

It would be interesting to know whether the same is true for more
general brackets, like Nambu-Jacobi brackets or brackets acting as
multidifferential operators. If we skip skew-symmetry, then it is
even not clear if the last ones have to be of first order. We can prove
the skew-symmetry for binary Nambu-Jacobi-Loday brackets and hope that the
methods used in the proof of an  algebraic version of the well-known
Kirillov's theorem on local Lie algebras (\cite{Gr}, Theorem 4.2) can
be of some help in proving a general result. We postpone these studies
to a separate paper.

\end{document}